\documentclass[preprint,12pt]{elsarticle}
\usepackage{lineno,hyperref}
\usepackage[utf8]{inputenc}
\usepackage[T1]{fontenc}    
\usepackage{lmodern}
\usepackage{amsfonts, latexsym, amsmath, amssymb}
\usepackage{mathtools}
\mathtoolsset{showonlyrefs}
\usepackage{graphicx}
\usepackage{microtype}
\usepackage{algorithm}
\usepackage{algpseudocode}
\usepackage{array}
\usepackage[figurename=Fig.]{caption}
\usepackage[margin=2cm]{geometry}
\biboptions{authoryear}

\makeatletter
\def\ps@pprintTitle{%
	\let\@oddhead\@empty
	\let\@evenhead\@empty
	\let\@oddfoot\@empty
	\let\@evenfoot\@oddfoot
}
\makeatother

\newcommand*{\horzbar}{\rule[.5ex]{3.0ex}{0.5pt}}

\begin{document}

\begin{frontmatter}
	
\title{A Reinforcement Learning Approach to the \\Stochastic Cutting Stock Problem}

\author[anselmoaddress,anselmoaddress2,opl]{Anselmo R. Pitombeira-Neto\corref{mycorrespondingauthor}}
\ead{anselmo.pitombeira@ufc.br}
\ead[url]{http://www.opl.ufc.br/authors/anselmo/}

\author[anselmoaddress2,opl]{Arthur H. F. Murta}

\address[anselmoaddress]{Department of Industrial Engineering, Federal University of Ceará, Campus do Pici, Fortaleza, Brazil}
\address[anselmoaddress2]{Graduate Program in Modeling and Quantitative Methods, Federal University of Ceará, Fortaleza, Brazil}
\address[opl]{OPL Lab - Operations Research in Production and Logistics, Federal University of Ceará, Fortaleza, Brazil}



\begin{abstract}
We propose a formulation of the stochastic cutting stock problem as a discounted infinite-horizon Markov decision process. At each decision epoch, given current inventory of items, an agent chooses in which patterns to cut objects in stock in anticipation of the unknown demand. An optimal solution corresponds to a policy that associates each state with a decision and minimizes the expected total cost. Since exact algorithms scale exponentially with the state-space dimension, we develop a heuristic solution approach based on reinforcement learning. We propose an approximate policy iteration algorithm in which we apply a linear model to approximate the action-value function of a policy. Policy evaluation is performed by solving the projected Bellman equation from a sample of state transitions, decisions and costs obtained by simulation. Due to the large decision space, policy improvement is performed via the cross-entropy method.  Computational experiments are carried out with the use of realistic data to illustrate the application of the algorithm. Heuristic policies obtained with polynomial  and Fourier basis functions are compared with myopic and random policies. Results indicate the possibility of obtaining policies capable of adequately controlling inventories with an average cost up to 80\% lower than the cost obtained by a myopic policy.
\end{abstract}

\end{frontmatter}


\section{Introduction}
The cutting stock problem (CSP) is a well studied problem in combinatorial optimization with applications in many industries, such as: textile, furniture, paper, glass, construction and manufacturing. In its basic setting, the problem consists in cutting smaller items from larger objects in stock, given the demand from customers, with the objective of minimizing trim loss. \citep{Gilmore1,Gilmore2} The problem has many variants which reflect specific features of real applications. \citep[see e.g.][]{Wascher,LEUNG,DELVALLE, FERNANDEZ}

The classic mathematical formulation of the CSP is the following: Let $d_{i}$ be the demand for an item of type $i \in \{1, 2, \dots, m\}$. The items are obtained by cutting larger objects in stock according to predefined cutting patterns. A cutting pattern $j \in \{1, 2, \dots, n\}$ is defined by a vector $a_j = (a_{1j}, a_{2j}, \dots, a_{mj})$ in which $a_{ij}$ is the number of items of type $i$ produced when a stock object is cut according to pattern $j$. The problem is then to determine how many objects $x_j$ to cut according to pattern $j$ in order to meet the demand while minimizing total trim loss, and can be formulated as the following integer linear program:
\begin{align}
\text{min} \quad & \sum_{j=1}^n g_j x_j, \label{eq:CSP_1}\\
 \text{s.t.} \quad & \notag \\
 & \sum_{j=1}^n a_{ij} x_j \geq d_i, \quad i \in \{1,2,\dots,m\}, \label{eq:CSP_2}\\
 & x_j \in \mathbb{Z}_+, \label{eq:CSP_3}
\end{align}
in which $g_j$ is the trim loss cost associated with pattern $j$. The CSP can be solved to optimality by using branch-and-bound algorithms.

Although mathematically elegant, the classic formulation of the CSP given in Eqs. \eqref{eq:CSP_1} - \eqref{eq:CSP_3} has limited applicability in industry, since it ignores key characteristics which occur in practice. First, it assumes that the demand is known at the time of decision, but in reality it is often unknown. Second, it is a static decision problem, while cutting stock operations in industry are typically dynamic, in which companies have to decide on a daily or weekly basis how many items to cut. Third, companies frequently keep items in inventory, so that, at a given decision time, there is an initial inventory of items that must be taken into account in the decision of cutting new items to fulfill the demand. This means that, in practice, cutting stock problems are mixed with lot-sizing problems with the added complication that production of individual items is coupled by the cutting patterns. We call this decision problem under a dynamic and uncertain scenario the \emph{stochastic cutting stock problem} (stochastic CSP).

We approach the stochastic CSP from a stochastic optimal control standpoint. This is, we seek a decision policy which specifies at each decision epoch how many objects in stock to cut in each pattern, given current inventory levels of demanded items. We present an unprecedented formulation of the stochastic CSP as a discounted infinite-horizon Markov decision process. Ideally, we would like to find an optimal policy which minimizes expected discounted total cost. Although in theory we may obtain an optimal policy through exact methods, in practice it will likely be computationally infeasible to obtain a guaranteed optimal policy due to the large size of the state and decision spaces. To overcome this difficulty, we propose an approximate solution based on reinforcement learning techniques.

In particular, we develop a model-free off-policy approximate policy iteration algorithm. The action-value function is approximated by a linear combination of basis functions. Samples of state transitions, decisions and rewards (costs in our case) are collected by simulation, and policy evaluation is performed by approximately solving the projected Bellman equation. Policy improvement is carried out by approximately solving a nonconvex nonlinear integer programming problem via the cross-entropy method. We investigate the proposed approach by means of computational experiments with the use of realistic data.

The structure of this paper is the following: In Section 2, we comment on related work; in Section 3, we formulate the stochastic CSP as a discounted infinite-horizon Markov decision process; in Section 4, we detail the proposed approximate policy iteration solution approach; in Section 5, we discuss computational results; finally, in Section 6 we draw some concluding remarks.

\section{Related Work}
Reinforcement learning (RL) is a computational approach to sequential decision problems based on the theory of Markov decisions processes and tries to overcome the well known \emph{curse of dimensionality} through approximations of value functions or policies \citep{puterman,suttonBarto}. RL has recently shown impressive results in video games \citep{mnih2015human}, board games \citep{silverGo} and robotics \citep{RL_Robotica}. It is closely related to approximate dynamic programming (ADP), a field with similar theory and methods but which evolved with different target applications \citep{powell_book, bertsekasNeuro}.

There has been increasing interest in applying RL and ADP methods to stochastic combinatorial optimization problems which arise in industry. The applications have spanned diverse areas, such as: transportation \citep{powell_vehicles}, scheduling \citep{lopes, jobshop}, energy storage management \citep{schneider, powell_ieee}, ambulance dispatching \citep{maxwell}, supply chain management \citep{kara, pourmoyaed} and blood bank management \citep{abdulwahab}.

However, as far as we know, RL has not been applied to the stochastic CSP. Some work has addressed the CSP in a dynamic but deterministic environment, in which demand for items is known in advance and supplied in a make-to-order manner. Since the demand may often not be supplied in a single time period, cutting operations have to be scheduled so as to optimize makespan, tardiness or other performance measure. In this case, it has been referred to in the literature as a \emph{combined cutting stock and scheduling problem} \citep{trkman2007, reinertsen, arbib, pitombeira2015, Pitombeira2019}.  Other papers address the CSP in a dynamic deterministic make-to-stock setting, in which point forecasts of future demands are available for a finite horizon and decisions on how manye objects to cut are taken in advance. \citep[See e.g.][]{nonas2000, poldi, melega, gramani2006,durak}

Just a few works approach a stochastic formulation of the CSP. The earliest we found is due to \cite{sculli}, who consider the sizes of the objects to be cut as random variables with known probability distributions. The problem is then to determine the position of the first cut of a fixed set of knives so as to minimize expected trim loss. \cite{krichagina} approach a CSP in the paper industry in which production and cutting decisions are 
made in a make-to-stock fashion. In a first stage, the average frequencies with which the cutting patterns are used are determined so as to meet average demand. Then, in a second stage, a scheduling policy decides if it shuts down a paper cutting machine to avoid building up too much inventory.

Some works apply the stochastic programming framework, which is fundamentally different from the stochastic control approach. \cite{alem} formulate a static CSP in which demands for items are unknown as a two-stage stochastic program with recourse. In the first stage, before seeing the demand, decisions are made on the items to be cut. In the second stage, after seeing the demand, costs are incurred for holding inventory or a penalty is incurred if demand is not met. The objective is to minimize the total expected cost in both stages. \cite{beraldi} propose a similar two-stage stochastic program with the difference that only the patterns to be used are decided upon before seeing the demand, while the actual amounts of objects cut according to each pattern are decided only after seeing the demand. \cite{zanarini} also proposes a two-stage stochastic programming model in which demands are uncertain. In the first stage, decisions are made on the sizes of the objects to be used in the cutting process, while in the second stage, after demands are known, the classic deterministic CSP is solved. It is worth noting that these three works consider only static one-time decisions.

In Section \ref{sec:markov}, we propose an unprecedented formulation of the stochastic CSP as a discounted infinite-horizon  Markov decision process. The main novelty of this formulation is that, in contrast to previous works, it takes into account both the dynamic and stochastic nature of the problem in practice.

\section{A Markov Decision Process Formulation}
\label{sec:markov}
Consider a cutting stock problem in a dynamic stochastic environment in which there are $m$ different items that can be demanded and cut from larger objects in stock, and there are $n$ different cutting patterns, where $a_{ij} $ is the number of items of type $i \in \{1,2, \dots, m\} $ obtained by cutting an object in the pattern $j \in \{1,2, \dots, n\}$. At each decision epoch $t$, there is an initial inventory $s_{it}$ for each type of item and a decision $x_{jt}$ must be made on how many objects to cut in each pattern before the demand is known, constrained by the availability of objects.

After the decision is taken, the amount of items produced is added to the initial inventory making up the available inventory to meet the demand. The demanded quantities deplete the available inventory, which results in the final inventory which will also be the initial inventory at subsequent time $t+1$. Demand which is not fulfilled is lost. There are costs related to trim loss, holding inventory to the next time period and lost sales. We would like to make decisions on which item quantities to cut at each decision epoch so as to minimize expected total cost over time. We name this problem the \emph{stochastic cutting stock problem}, since the demand and costs are not fully known at the decision time.

We formulate the stochastic CSP as a discounted infinite-horizon Markov decision process (MDP). Let $s_t \in \mathcal{S} $ be a vector corresponding to the state of the system at time $t$ and $\mathcal{S}$ a finite set of states. We define $s_t = (s_ {1t}, s_ {2t}, \dots, s_ {mt})$ where $s_{it} \in \{0, \dots, s_{\mathrm {max}} \}$ as the initial inventory of items of type $i$ at time $ t $. Let $x_t = (x_{1t}, x_{2t}, \dots, x_ {nt})$ be a vector corresponding to the decision (we also use the term \emph{action} interchangeably) at time $t$ when the state is $s_t = s$, in which $x_{jt} \in \{0, 1, \dots, x_{\mathrm{max}} \}$ corresponds to the number of objects in stock that are cut in  pattern $j$. Decisions are constrained to assume values in a set $\mathcal {X}_s$ of feasible decisions given the state $s$ and we also define the set of all possible decisions as $\mathcal{X} = \bigcup_{s \in \mathcal{S}} \mathcal{X}_s$. We denote as $d_{t + 1} = (d_{1, t + 1}, d_{2, t + 1}, \dots, d_{m, t + 1})$ the demand vector for items that occurs in the interval $ (t, t + 1] $. We assume that this demand becomes known only \emph {after} the decision to cut the objects in stock, so that $ d_{t + 1} $ is a random discrete vector with a conditional probability function $\mathbb{P}(d_{t + 1} | s_t, x_t)$.

The state transition function is defined by $s_{t+1} = f(s_t,x_t,d_{t+1})$, in which $s_{t+1}$ is the initial inventory at time $t+1$, whose components $s_{i,t+1}$  are given by
\begin{equation}
    s_{i,t+1} = \Bigg[s_{it}+\sum_{j = 1}^n a_{ij}x_{jt}-d_{i,t+1}\Bigg]^+, \quad i \in \{1, \dots, m\}, \label{eq:funcao_transicao}
\end{equation}
in which $[x]^+ = \max(0,x)$. Notice in \eqref{eq:funcao_transicao} that
\begin{equation}
\sum_{j = 1}^n a_{ij}x_{jt} \notag
\end{equation}
corresponds to the quantity of item $i$ added to the inventory by cutting $x_{jt}$ stock objects according to patterns $j \in \{1,2,\dots, n\}$ at time $t$. Notice also that $s_{t+1}$ is a random vector, since it is a function of the random vector $d_{t+1}$.

The cost function (also called \emph{reward function} in reinforcement learning) is the sum of costs associated with trim loss, holding inventory  and lost sales, given by
\begin{equation}
    c(s_t,x_t,d_{t+1}) = \sum_{j=1}^n g_{j}x_{jt}+\sum_{i=1}^m h_{i}^+\Bigg[s_{it}+\sum_{j = 1}^n a_{ij}x_{jt}-d_{i,t+1}\Bigg]^+ +\sum_{i=1}^m h_{i}^-\Bigg[d_{i,t+1} -\Bigg (s_{it}+\sum_{j = 1}^n a_{ij}x_{jt}\Bigg)\Bigg]^+, \label{eq:funcao_custo}
\end{equation}
in which $g_j$ is a cost term associated with trim loss by using pattern $j$, $h_i^+$ is the inventory holding cost of item $i$ for a time period, and $h_i^-$ is the lost sales cost for item $i$. Notice that the third term in the right-hand side of \eqref{eq:funcao_custo} corresponds to the total cost of lost sales, which is incurred when
\begin{equation}
d_{i,t+1} -\Bigg (s_{it}+\sum_{j = 1}^n a_{ij}x_j\Bigg) > 0, \notag
\end{equation}
i.e., when the demand for item $i$ is greater than the available inventory. A noteworthy characteristic of the cost function \eqref{eq:funcao_custo} is that it is random at the decision time $t$, since the decision maker (also called the agent in reinforcement learning) does not know the demand $d_{t+1}$ which will be realized during time interval $[t,t+1)$.

Let $\pi$ be a stationary Markovian policy defined by $\pi := \{x(s)\}$, in which $x: \mathcal{S} \to \mathcal{X}$ is a function that associates with each state $s \in \mathcal{S}$ a decision $x = x(s)$. (Notice that we can drop the subscript $t$ by convenience since we henceforth assume stationarity of the system and the policy.) A feasible policy associates with each state $s \in \mathcal{S}$ a decision $x$ in the set $\mathcal{X}_s$ of feasible decisions given by
\begin{equation}
\mathcal{X}_s = \left\{
\begin{aligned}
    & s_i+\sum_{j=1}^n a_{ij}x_j  \leq s_{\mathrm{max}}, \quad i=1,2, \dots, m,\\
    & \sum_{j=1}^n x_j \leq x_\mathrm{max},\\
    & x_j \in \mathbb{Z}_+, \quad j =1, 2, \dots, n
\end{aligned}
\right\},    \label{eq:MDP_constraints}
\end{equation}
in which $s_\text{max}$ is the maximum inventory and $x_\text{max}$ is the maximum number of objects in stock which may be cut in a given time period.

Given an initial state $s_0 = s$, we define as $v_\pi: \mathcal{S} \to \mathbb{R}$ a value function associated with a policy $\pi$:
\begin{equation}
    v_\pi(s) := \lim_{T \to \infty} \mathbb{E} \Bigg[\sum_{t=0}^T \gamma^t c(s_t,x(s_t),d_{t+1})|s_0 = s \Bigg], \quad \forall s \in \mathcal{S}, \label{eq:funcao_valor}
\end{equation}
in which $0<\gamma < 1$ is a \emph{discount factor} and the expected value is computed with regards to the joint probability function of the states $(s_1, s_2, ...\dots, s_T)$ induced by policy $\pi$. The \emph{Markov decision problem} corresponds to finding an optimal policy $\pi^\star$ which minimizes the value function \eqref{eq:funcao_valor} for all states $s \in \mathcal{S}$:
\begin{equation}
    \pi^\star \in \arg \min_{\pi \in \Pi_{\text{md}}} v_\pi(s), \quad \forall s \in \mathcal{S}, \label{eq:funcao_valor_2}
\end{equation}
in which $\Pi_{\text{md}}$ denotes the class of deterministic Markovian policies. Notice that the problem \eqref{eq:funcao_valor_2} cannot in general be solved directly. However, value functions satisfy the well known Bellman equation, which in our case is given by
\begin{equation}
    v_\pi(s) =   \mathbb{E}[c(s, x(s),d)+ \gamma v_{\pi}(s') |s], \quad \forall s \in \mathcal{S}, \label{eq:bellman}
\end{equation}
in which the expected value in \eqref{eq:bellman} is computed with regards to the conditional probability function $\mathbb{P}(d|s,x)$ of the demand $d$ and $s' = f(s,x,d)$ is given by Eq. \eqref{eq:funcao_transicao}.

In theory, an optimal policy may be obtained by exact methods such as value iteration, policy iteration or linear programming \citep{puterman}. We are particular interested in the policy iteration method \citep{howard},  which swaps between two steps: policy evaluation and policy improvement.  Given an arbitrary policy $\pi^{(k)}$ at iteration $k$ of the method, in the evaluation step we obtain its value function $v_{\pi^{(k)}}$ by solving Bellman equation \eqref{eq:bellman}. Then in the improvement step, a new policy $\pi^{(k+1)} := \{x^{(k+1)}(s)\}$ is obtained by acting greedily with respect to the value function $v_{\pi^{(k)}}$:
\begin{equation}
    x^{(k+1)}(s) \in \arg \min_{x \in \mathcal{X}_s} \mathbb{E}[c(s, x,d)+ \gamma v_{{\pi}^{(k)}}(s') |s], \quad \forall s \in \mathcal{S}.    \label{eq:politica_gulosa} \notag
\end{equation}
It can be shown \citep{puterman} that $\pi^{(k+1)}$ is better than $\pi^{(k)}$ in the sense that $v_{\pi^{(k+1)}}(s) \leq v_{\pi^{(k)}}(s), \quad \forall s \in \mathcal{S}$. ($\geq$ if maximizing rewards.) Policy iteration begins with an initial policy $\pi^{(0)}$ and generates a sequence of policies $\pi^{(0)}, \pi^{(1)}, \pi^{(2)}, \dots$ which stops when $\pi^{(k+1)} = \pi^{(k)}$. Value iteration may be regarded as a particular case of policy iteration in which the policy evaluation and improvement steps are merged in a single step.

Although policy iteration is computationally more efficient than alternative exact methods for solving MDPs, such as linear programming \citep{suttonBarto}, it is impractical for problems with large states spaces. Its main limitation is that the value function must be computed for all states, and in vector-valued state variables the number of states increases exponentially with the dimension of the vector. (The well known \emph{curse of dimensionality}.) Other limitation of policy iteration is that it relies on the exact computation of the expectation in the Bellman equation \eqref{eq:bellman}, which will not be possible in general. 

Since exact policy iteration will be computationally infeasible in our case, in the next section we leverage methods developed in the field of reinforcement learning to obtain an approximate solution to the stochastic CSP.

\section{An Approximate Policy Iteration Solution Approach}
In this section, we develop an approximate policy iteration method to solve the stochastic CSP. Exact policy iteration alternates between two steps: policy evaluation and policy improvement. However, in approximate policy iteration these two steps are intertwined and the boundaries of each step are not so crispy.

In order to cope with the large state space, we use a parameterized value-function approximation defined by a linear combination of basis functions. Given a greedy policy with respect to a current value function, we gather samples from state transitions, decisions and costs by simulation. Computing the greedy action is nontrivial, since the decision space is also very large and corresponds to solving a nonconvex nonlinear integer programming problem. We then solve it heuristically by using the cross-entropy method. This sampling step corresponds to the policy improvement step.

\begin{figure}
	\centering
	\includegraphics[scale=0.70, trim=0 19cm 0 2cm 0,clip]{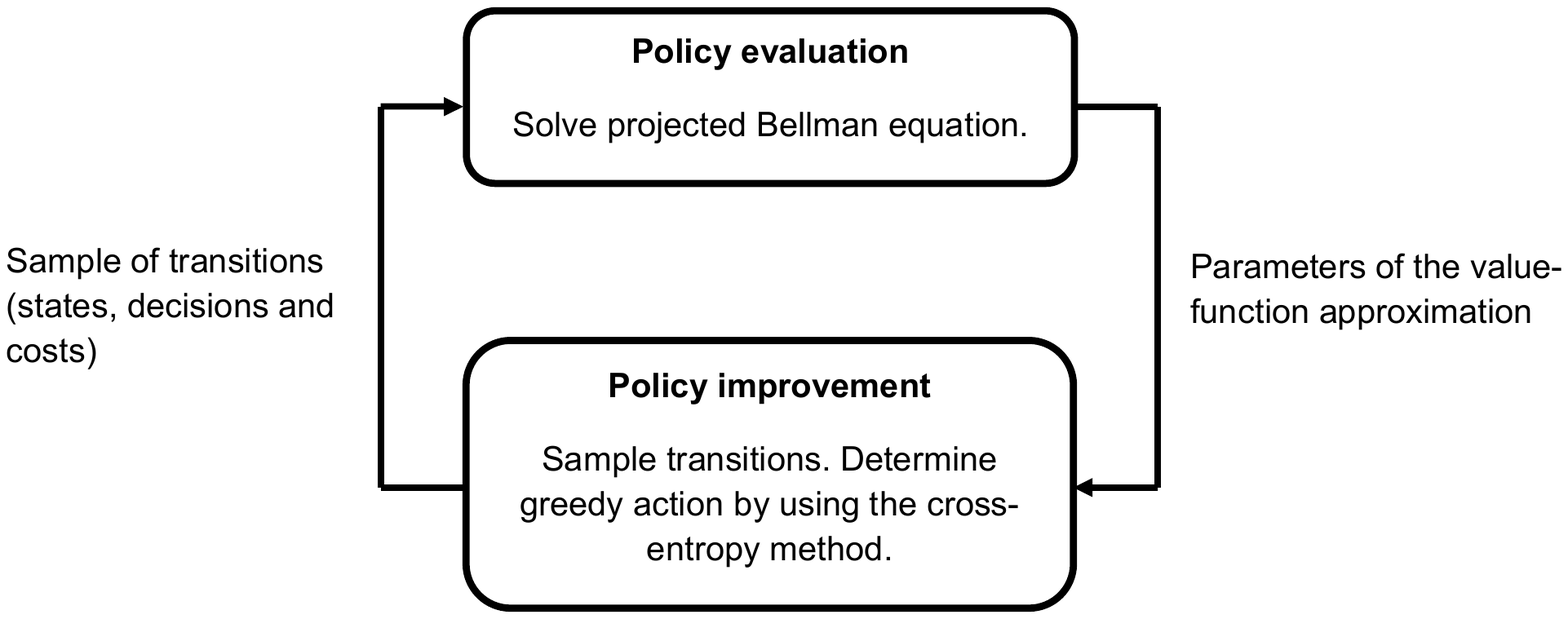}
	\caption{General scheme of the proposed approximate policy iteration algorithm for the stochastic cutting stock problem}
	\label{fig:api}
\end{figure}

Given the collected sample of transitions, the policy evaluation step is accomplished to solving the Bellman equation. As we use a value-function approximation, we solve the projected Bellman equation, whose solution is the fixed point of the Bellman operator projected into the space spanned by the basis functions. Fig. \ref{fig:api} illustrates the general sketch of the proposed algorithm. We detail its development in the next sections.
\subsection{Policy evaluation via projected Bellman equation}
\label{sec:policy_evaluation}
Instead of working with \emph{state-value functions} $v_\pi(s)$, we will work with \emph{action-value functions} (also known as q-functions), defined as
\begin{equation}
    q_\pi(s,x) := \lim_{T \to \infty} \mathbb{E} \Bigg[\sum_{t=0}^T \gamma^t c(s_t,x(s_t),d_{t+1})|s_0 = s, x(s_0) = x \Bigg], \quad \forall s \in \mathcal{S}, \label{eq:funcao_valor_acao} \notag
\end{equation}
which corresponds to the value of starting at state $s$, choosing initial action $x$ and then following policy $\pi := \{x(s)\}$ for $t \in \{1, 2, \dots\}$.
Action-value functions satisfy a form of Bellman equation, given by
\begin{equation}
    q_\pi(s,x) = \mathbb{E}[c(s, x,d)+ \gamma q_{\pi}(s',x(s')) |s,x], \quad \forall s \in \mathcal{S}, \forall x \in \mathcal{X}, \label{eq:bellman_q} \notag
\end{equation}
in which the expected value is computed with regards to the conditional probability function $\mathbb{P}(d|s,x)$ and $s' = f(s,x,d)$ is the transition function. The use of action-value functions simplifies the computation of a greedy policy, since in this case a greedy action may be determined without knowing the MDP model or computing an expectation. Given a state $s \in \mathcal{S}$, a greedy action $x(s)$ may be determined by
\begin{equation}
    x(s) \in \arg \min_{x \in \mathcal{X}_s} q_\pi(s,x).     \label{eq:politica_gulosa_q}
\end{equation}
In order to cope with the large state space of the stochastic CSP, we use a linear model to approximate the action-value function of a policy $\pi$:
\begin{equation}
    q_\pi(s,x;\theta) = \sum_{k = 1}^K \phi_k(s,x)\theta_k, \label{eq:modelo_linear}
\end{equation}
in which $\phi_k(s,x)$ are $K$ basis functions, which represent \emph{features} associated with state $s$ and action $x$, and $\theta = (\theta_1, \theta_2, \dots, \theta_K)$ is a parameter vector.

Our objective is then to \emph{train} the model \eqref{eq:modelo_linear} (this is, to estimate its parameters) from samples of state transitions, actions and costs obtained by simulation. Notice though that the exact action-value function does not necessarily lie in the space spanned by the basis functions. Proposed criteria to set the parameters include Bellman error minimization and projected Bellman error minimization \citep{lagoudakis2003least}. We adopt the projected Bellman error criterion to choose suitable values for the parameters.

Consider the \emph{Bellman operator} $T_\pi : \mathcal{Q} \to \mathcal{Q}$ associated with a policy $\pi :=\{x(s)\}$, in which  $\mathcal{Q}$ is the space of action-value functions, defined by
\begin{align}
    (T_\pi q)(s,x) & := \mathbb{E}[c(s, x,d)+ \gamma q(s',x(s')) |s,x], \quad \forall s \in \mathcal{S}, \forall x \in \mathcal{X}, \notag \\
    & := \sum_{d \in \mathcal{D}} \mathbb{P}(d|s,x) (c(s,x,d)+ \gamma q(s',x(s'))), \label{eq:Bellman_operator}
\end{align}
in which $\mathcal{D}$ is the support of the probability function $\mathbb{P}(d|s,x)$. Notice that $q$ is an arbitrary action-value function in the space $\mathcal{Q}$. It will be convenient to write the Bellman operator in vector notation:
\begin{equation}
    T_\pi(q) := \bar{c}+\gamma \bar{q}, \notag
\end{equation}
in which $\bar{c} \in \mathbb{R}^{|\mathcal{S} \times \mathcal{X}|}$ is a column-vector, whose components are given by
\begin{equation}
 \bar{c}(s,x) = \sum_{d \in \mathcal{D}} \mathbb{P}(d|s,x) c(s,x,d),   \notag
\end{equation}
and $\bar{q} \in \mathbb{R}^{|\mathcal{S} \times \mathcal{X}|}$ is a column-vector, whose components are given by
\begin{equation}
 \bar{q}(s,x) = \sum_{d \in \mathcal{D}} \mathbb{P}(d|s,x) q(s', x(s')), \notag
\end{equation}
and $s' = f(s,x,d)$ is the transition function given by \eqref{eq:funcao_transicao}.

Now, given an arbitrary policy $\pi$, we want to obtain its action-value function $q_\pi$. We first notice that, from Bellman equation \eqref{eq:bellman_q}, the action-value function $q_\pi$ is the fixed point of the Bellman operator \eqref{eq:Bellman_operator}. However, we cannot solve directly for the fixed point, since we are using an approximation for the action-value function given by \eqref{eq:modelo_linear} and there is no guarantee that $q_\pi$ lies in the space spanned by the basis functions $\phi_k, k \in \{1,2,\dots, K\}$. Instead, we solve the \emph{projected Bellman equation} \citep{bertsekasNeuro}:
\begin{equation}
    \Phi \theta = \Pi_{\xi}T_\pi (\Phi \theta), \label{eq:bellman_projetada}
\end{equation}
in which $\Phi$ is a matrix with $|\mathcal{S}| \times \sum_{s \in \mathcal{S}}|\mathcal{X}_s|$ rows and $K$ columns whose rows are the feature column-vectors $\phi(s,x)^T = (\phi_1(s,x),\phi_2(s,x), \dots, \phi_K(s,x))$ for each pair $(s,x)$ with $s \in \mathcal{S}$ and $x \in \mathcal{X}_s$:
\begin{equation}
\Phi =
\left[
\begin{array}{c}
\horzbar  \, \phi(s_1,x_1)^T \,  \horzbar \\
\horzbar \,  \phi(s_1,x_2)^T \,  \horzbar\\
\vdots \\
\horzbar \,  \phi(s_{|\mathcal{S}|},x_{|\mathcal{X}_s|})^T \,  \horzbar
\end{array}\right],
\end{equation}
and $\Pi_{\xi}$ is a projection matrix to the space spanned by the basis functions in which the projection is defined with relation to a weighted norm   $||.||_\xi$ and $\xi$ is the weight vector. In the particular case of a Euclidean norm, the projection matrix corresponds to the least squares solution, given by \citep{geramifard}:
\begin{equation}
 \Pi_{\xi} = \Phi (\Phi^T\Xi\Phi)^{-1}\Phi^T\Xi, \label{eq:matriz_projecao}
\end{equation}
in which $\Xi = \mathrm{diag}(\xi)$.
Substituting \eqref{eq:matriz_projecao} in \eqref{eq:bellman_projetada}, we have
\begin{align}
\Phi \theta & = \Phi (\Phi^T\Xi\Phi)^{-1}\Phi^T\Xi(\bar{c} +\gamma \bar{\Phi} \theta), \notag \\
\Phi^T\Xi\Phi \theta & = \Phi^T\Xi(\bar{c} +\gamma \bar{\Phi} \theta), \notag\\
\Phi^T\Xi(\Phi - \gamma \bar{\Phi})\theta & = \Phi^T\Xi\bar{c}, \label{eq:erro_bellman_projetado}
\end{align}
in which $\bar{\Phi}$ is a matrix with $|\mathcal{S}| \times \sum_{s \in \mathcal{S}}|\mathcal{X}_s|$ rows and $K$ columns whose rows $\bar{\phi}(s,x)$ are given by
\begin{equation}
    \bar{\phi}(s,x) = \sum_{d \in \mathcal{D}} \mathbb{P}(d|s,x) \phi(s', x(s')), \quad s' = f(s,x,d), \notag
\end{equation}
so that
\begin{equation}
\bar{\Phi} =
\left[
\begin{array}{c}
\horzbar  \, \bar{\phi}(s_1,x_1)^T \,  \horzbar \\
\horzbar \,  \bar{\phi}(s_1,x_2)^T \,  \horzbar\\
\vdots \\
\horzbar \,  \bar{\phi}(s_{|\mathcal{S}|},x_{|\mathcal{X}_s|})^T \,  \horzbar
\end{array}\right]. \notag
\end{equation}
Notice that \eqref{eq:erro_bellman_projetado} corresponds to solving a linear system of equations $A\theta = b$ in which $A = \Phi^T\Xi(\Phi - \gamma \bar{\Phi})$ is a $K \times K$ matrix and $b = \Phi^T\Xi\bar{c}$ is a $K \times 1$ column-vector. Notice also that the linear system \eqref{eq:erro_bellman_projetado} depends on all states $s \in \mathcal{S}$ and actions $x \in \mathcal{X}$. However, we can solve it approximately if we have a sample of $N$ transitions $<s_t,x_t,c_{t+1},s_{t+1}>$, $t = 0, \dots, N-1$. We then form the approximate matrix $\hat{A}$ and vector $\hat{b}$:
\begin{align}
    \hat{A} & = \frac{1}{N} \sum_{t=0}^{N-1} \phi(s_t,x_t)(\phi(s_t,x_t) - \gamma \phi(s_{t+1}, x(s_{t+1})))^T, \notag \\
    \hat{b} & = \frac{1}{N} \sum_{t=0}^{N-1} \phi(s_t,x_t) c_{t+1}, \notag
\end{align}
in which $\phi(s,x)$ is the feature column-vector and transitions are sampled according to a probability distribution corresponding to the weights $\xi$. Therefore, we can obtain the parameter vector by solving the approximate linear system
\begin{equation}
\hat{A} \theta = \hat{b}. \label{eq:linear_system}
\end{equation}
Let $\hat{\theta}$ be a solution to the linear system \eqref{eq:linear_system}. (We may use the pseudo-inverse if $\hat{A}$ is singular.) Then the action-value function of policy $\pi$ is approximated by
\begin{equation}
    q_\pi(s,x;\hat{\theta}) = \sum_{k = 1}^K \phi_k(s,x)\hat{\theta}_k. \notag
\end{equation}
In summary, evaluating a policy $\pi$ reduces to solving the linear system \eqref{eq:linear_system} from a sample of transitions taken while using policy $\pi$. In the next section \ref{sec:policy_improvement}, we describe how to obtain an improved policy.

\subsection{Policy improvement via the cross-entropy method}
\label{sec:policy_improvement}
Given a parameter vector $\theta$, an improved policy $\pi'$ relative to a current policy $\pi$ may be obtained by acting greedily with respect to the approximate action-value function $q_\pi(s,x;\theta$):
\begin{equation}
    x(s) \in \arg \min_{x \in \mathcal{X}_s} \sum_{k = 1}^K \phi_k(s,x)\theta_k, \quad \forall s \in \mathcal{S}. \label{eq:politica_gulosa_q_aprox}
\end{equation}
Notice that, in the case of the stochastic CSP, \eqref{eq:politica_gulosa_q_aprox} will be in general a nonconvex integer nonlinear programming problem. These kind of mathematical programs are among the hardest to solve exactly. Currently available exact solvers are not sufficiently flexible or computationally efficient for our application, since we will have to solve problem \eqref{eq:politica_gulosa_q_aprox} thousands of times during the simulation of a sample of transitions to evaluate a single policy. Moreover, solving \eqref{eq:politica_gulosa_q_aprox} only approximately may not be detrimental to performance, since we are already working with an approximation of the action-value function.

We use the cross-entropy method to obtain a heuristic solution to \eqref{eq:politica_gulosa_q_aprox}. The cross-entropy method employs a parametric probability distribution over the space of solutions, which is used to generate a sample of candidate solutions. The objective values of the candidate solutions are computed and a fraction of the best solutions (the elite group) is selected. The elite group is then used to estimate new parameters to the probability distribution. This process is repeated until some stopping criterion is achieved. A detailed explanation of the method is given by \cite{cross_entropy}.

Given a current state $s$, we generate candidate feasible decisions in the following way: We first sample the total number $x_\text{total}$ of objects which will be cut from a discrete uniform distribution:
\begin{equation}
x_\text{total} \sim \mathtt{DiscUnif}(0,x_\text{max}), \notag
\end{equation}
in which $x_\text{max}$ is the maximum number of available objects in stock. We then use a multinomial probability distribution to generate a candidate solution $x$:
\begin{equation}
 x \sim \mathtt{Multinomial}(x_\text{max}; p_1, p_2, \dots, p_n), \label{eq:multinomial}
\end{equation}
in which $p_j$ is the probability of using a cutting pattern $j \in \{1,2,\dots,n\}$. Feasibility of the sampled candidate is then checked against the set $\mathcal{X}_s$. Infeasible candidate solutions are rejected.

Let $x^{(1)}, x^{(2)}, \dots, x^{(N)}$ be a sample of candidate solutions taken from \eqref{eq:multinomial}. We compute the action values $q_\pi(s,x^{(1)}; \theta), q_\pi(s,x^{(2)};\theta),\dots, q_\pi(s,x^{(N)};\theta)$ from \eqref{eq:modelo_linear}. We then sort the action values in increasing order. Let $\delta$ be the $\lceil \rho \times  N \rceil$-th order statistic, with $0 < \rho < 1$, such that the $\lceil \rho \times  N \rceil$ best candidate solutions have action values $q_\pi(s,x;\theta) \leq \delta$. (For example $\rho = 0.10$ select the $10\%$ best candidate solutions, also called elite solutions.) We then determine updated parameters $p' = (p'_1, p'_2, \dots, p'_n)$ of the multinomial distribution by maximizing the cross-entropy function, which in this case will be mathematically equivalent to determining the maximum likelihood estimate given only the elite solutions:
\begin{align}
\max \quad &  \frac{1}{N} \sum_{i=1}^{N} I\{q_{\pi}(s,x^{(i)};\theta) \leq \delta\} \ell(x^{(i)}; p), \label{eq:cross-entropy} \\
\text{s.t.} \quad & \notag \\
\quad &  \sum_{j=1}^n p_j = 1, \notag\\
\quad & p_j \geq 0 \quad j \in \{1, 2, \dots, n\}, \notag
\end{align}
in which $I\{.\}$ is an indicator function and $\ell(x^{(i)}; p)$ is the multinomial log-likelihood, given by
\begin{equation}
\ell(x; p) = \sum_{j=1}^n x_j \ln p_j. \notag
\end{equation}
It can be shown, through the method of Lagrange multipliers, that $p' = (p'_1, p'_2, \dots, p'_n)$ which maximizes \eqref{eq:cross-entropy} is given by
\begin{equation}
p_j' =  \frac{ \sum_{i = 1}^N I\{q_\pi(s,x^{(i)}; \theta) \leq \delta\} x^{(i)}_j}{\sum_{j=1}^n \sum_{i = 1}^N I\{q_\pi(s,x^{(i)};\theta) \leq \delta\} x^{(i)}_j}, \quad  j \in \{1,2,\dots,n\}. \label{eq:p_estimator}
\end{equation}
Although a bit daunting, Eq. \eqref{eq:p_estimator} means that $p'_j$ is given by the sample frequency with which cutting pattern $j$ was used by the elite candidate solutions. Initial probabilities are set at $p_j = 1/n, j \in \{1,2,\dots,n\}$. Algorithm \ref{alg:cross-entropy} summarizes the steps to obtain a heuristic solution to \eqref{eq:politica_gulosa_q_aprox} with the use of the cross-entropy method.
\begin{algorithm}[t]
	\caption{Cross-entropy method for obtaining a heuristic greedy action}
	\label{alg:cross-entropy}
	\begin{algorithmic}[1]
		\State \textbf{input:} Current state $s$ and $x_\text{max}$, basis functions $\phi(.)$,\\ \qquad parameters $\theta$,  algorithm parameters  $N_1,N_2, \rho$;
		\State \textbf{initial step:} Set $p^{(0)}_j = 1/n, j \in \{1,2,\dots, n\}$, $q_\text{best} \gets +\infty$;
		\For{$k \gets 0 \dots N_1-1$}
		    \For{$i \gets 1 \dots N_2$}
		        \While{\text{True}}
		            \State Sample $x^{(k,i)}_\text{total} \sim \mathtt{DiscUnif}(0,x_\text{max})$;
		            \State Sample candidate $x^{(k,i)} \sim \mathtt{Multinomial}(x^{(k,i)}_\text{total};p^{(k)}_1, p^{(k)}_2, \dots,  p^{(k)}_n)$;
              		\If{$x^{(k,i)} \in \mathcal{X}_s$}
                        \State \textbf{break while};
		            \EndIf
		        \EndWhile
		        \If{$q_\pi(s, x^{(k,i)};\theta) < q_\text{best}$}  \Comment{Compute $q_\pi(s, x^{(k,i)};\theta)$ from equation \eqref{eq:modelo_linear}}
		        \State $x_\text{best} \gets x^{(k,i)}$; \Comment{Candidate solution turns to current best solution}
		        \State $q_\text{best} \gets q_\pi(s, x^{(k,i)};\theta)$;
		        \EndIf
		    \EndFor
		    \State Sort action values $q_\pi(s, x^{(k,i)};\theta)$, $i \in \{1,2,\dots, N_2\}$ in increasing order;
		    \State Let $\delta^{(k)}$ be the $\lceil \rho \times  N_2 \rceil$-th order statistic in the sorted list of action values;
		    \State Update probabilities
		    \[p^{(k+1)}_j =  \frac{ \sum_{i = 1}^{N_2} I\{q_\pi(s,x^{(k,i)}) \leq \delta^{(k)}\} x^{(k,i)}_j}{\sum_{j=1}^n \sum_{i = 1}^{N_2} I\{q_\pi(s,x^{(k,i)}) \leq \delta^{(k)}\} x^{(k,i)}_j}, \quad  j \in \{1,2,\dots,n\};\]
        \EndFor	
		\State \textbf{return} $x_\text{best}$ \Comment{Best found heuristic greedy action}
	\end{algorithmic}
\end{algorithm}

\subsection{Appproximate policy iteration algorithm}

Algorithm \ref{alg:iteracao_politica_aproximada} describes the proposed approximate policy iteration to obtain a heuristic policy to the stochastic CSP. The algorithm has two loops: in the outer loop, each iteration $i$ corresponds to the evaluation step of a greedy policy $\pi^{(i)}$ relative to the approximate action-value function given by the linear model \eqref{eq:modelo_linear}; in the inner loop, a sample of transition is taken by simulating states, decisions and costs. The sample is used to compute the matrix $\hat{A}$ and the vector $\hat{b}$ in order to obtain a new parameter vector $\theta^{(i+1)}$. The outer loop is run for $L_1$ policy iterations, while the inner loop is run for $L_2$ iterations, which correspond to the size of the sample of transitions. It returns all computed parameter vectors $\{\theta^{(1)},\theta^{(2)}, \dots, \theta^{(L_1)}\}$. This allows us to reevaluate all generated policies.The algorithm is \emph{off-policy}, i.e., state transitions can be sampled arbitrarily without following the stationary probability distribution induced by the current policy being evaluated. 
\begin{algorithm}[h!]
	\caption{Approximate policy iteration for the stochastic CSP}
	\label{alg:iteracao_politica_aproximada}
	\begin{algorithmic}[1]
		\Statex \textbf{input}: Basis functions $\phi(.)$, $\theta^{(0)}, \gamma, L_1, L_2$;
		\For{$i \gets 0 \dots L_1-1$}
		\Statex \quad \; \textbf{Simulation of the greedy policy relative to} $\theta^{(i)}$: \Comment{Policy improvement}
		\State $\hat{A}^{(i,0)} \gets 0, \hat{b}^{(i,0)} \gets 0$;
		\For{$t \gets 0 \dots L_2-1$}
		\State Sample initial inventory $s_t \in \mathcal{S}$;
		\State Sample decision $x_t \in \mathcal{X}_{s_t}$;
		\State Sample demand $d_{t+1} \sim p(d|s_t,x_t);$
		\State Compute cost $c_{t+1} = c(s_t,x_t,d_{t+1})$ from equation  \eqref{eq:funcao_custo};
		\State Compute transition $s_{t+1} = f(s_t,x_t,d_{t+1})$ from equation \eqref{eq:funcao_transicao};
		\State Compute decision $x_{t+1} = x(s_{t+1})$, by solving
		\[
		x(s_{t+1}) \in \arg \min_{x \in \mathcal{X}_{s_{t+1}}} \sum_{k = 1}^K \phi_k(s_{t+1},x)\theta^{(i)}_k,
		\]
		\qquad \quad by means of the cross-entropy method (Algorithm \ref{alg:cross-entropy});
		\State Update matrices
		\begin{align}
		\hat{A}^{(i,t+1)} & \gets \hat{A}^{(i,t)}+ \phi(s_t,x_t)(\phi(s_t,x_t) - \gamma \phi(s_{t+1}, x_{t+1}))^T, \notag \\
		\hat{b}^{(i,t+1)}& \gets \hat{b}^{(i,t)}+\phi(s_t,x_t) c_{t+1}; \notag
		\end{align}
		\EndFor
		\Statex \quad \; \textbf{Solve projected Bellman equation:} \Comment{Policy evaluation}
		\State Determine $\theta^{(i+1)}$ by solving $\hat{A}^{(i,L_2)} \theta = \hat{b}^{(i,L_2)}$; 
		\Statex \quad \; (Use the pseudo-inverse if $\hat{A}^{(i,L_2)}$ is singular.)
		\State 
		\EndFor
		\State \textbf{return} $\{\theta^{(1)},\theta^{(2)}, \dots, \theta^{(L_1)}\}$
	\end{algorithmic}
\end{algorithm}

\subsection{Features and basis functions}
For the application of Algorithm \ref{alg:iteracao_politica_aproximada}, it is necessary to choose basis functions used in approximating the action-value function of a policy. There are many types of basis functions, but in our application we will focus on two particular ones: polynomial and Fourier basis functions, which we found to have convenient properties in our case. First, they are easy to specify and, second, they are computationally cheap to compute. This is critical, since we cannot spend too much time to obtain a greedy action when approximately solving \eqref{eq:politica_gulosa_q_aprox}. Polynomial basis functions are simply polynomial terms in the state $s$ and decision $x$, while Fourier basis functions are sinusoidal functions such as sines and cosines \citep{konidaris}.

In addition, instead of defining the features $\phi(s,x)$ as explicit functions of $s$ and $x$, we take an intermediate step and aggregate state and decision in a \emph{post-decision} state $s^x = f(s,x)$, in a way that each feature will be defined in terms of $s^x = (s_1^x, s_2^x, ..., s_m^x)$. In the case of the stochastic CSP, we define as the post-decision state $s_i^x$ related to item $i$ the available inventory, given by the sum of the initial inventory $s_i$ and the quantity of the item resulting from the decision $x$:
\begin{equation}
s_i^x = s_i +\sum_{j = 1}^n a_{ij} x_j. \label{eq:post-decision}
\end{equation}
The main advantage of working with the notion of the post-decision state is that it naturally takes into account the relationship between the inventory of items and the cutting patterns, besides being a more compact form which is independent of the number of cutting patterns.

In the case of a polynomial basis, the action-value function may be written as \citep{suttonBarto}:
\begin{equation}
q_{\pi}^\text{poly}(s,x;\theta) = \sum_{k=1}^K \prod_{i=1}^{m} (s_i^x)^{c_{ki}}\theta_k, \notag
\end{equation}
in which $c_{ki}$ are integers which specify the order of the polynomial features.

Regarding the Fourier basis, it uses as basis functions the terms of the Fourier series, which is commonly used to approximate continuous functions. In our case, we will use a Fourier basis with cosine functions only, which according to \cite{konidaris} is sufficient to approximate any function defined in a non-negative domain, given by
\begin{equation}
q_{\pi}^\text{Fourier}(s,x;\theta) = \sum_{k=1}^K \cos(\pi c_k \cdot s^\text{norm})\theta_k, \notag
\end{equation}
in which $c_k = (c_1, c_2, \dots, c_m)$ is a vector of integers which specify the frequencies of the cosine functions. Larger integer numbers generate cosine functions with higher frequencies. In addition, we use normalized post-decision states, such that $s^\text{norm} = (s_1^\mathrm{norm}, s_2^\mathrm{norm}, \dots, s_m^\mathrm{norm})$ and $s_i^\mathrm{norm} = s_i^x/s_\mathrm{max}$.

\section{Results and Discussion}
In this section, we investigate the application of Algorithm \ref{alg:iteracao_politica_aproximada} with the use of a realistic dataset. Our objective is to evaluate whether our approximate policy iteration algorithm can generate a decision policy for the stochastic CSP with acceptable performance, which we mean as a policy that is able of meeting the demand without maintaining excessive inventory levels. In other words, our \emph{agent} has to learn which cutting patterns to use over time so that it maintains sufficient and not very high inventory levels. We compare the obtained heuristic policy with myopic and random policies.

\subsection{A myopic policy}
We benchmark the trained policies with a myopic policy which is a reasonable heuristic to the stochastic CSP in practice. The proposed myopic policy decides which items to cut at each time period according to the expected demand. Then, at each time period, the myopic policy solves a deterministic CSP model so that the sum of the current initial inventory at hand and the items cut is greater than the expected demand at minimal trim loss. Let $\bar{d} = (\bar{d}_1, \bar{d}_2, \dots, \bar{d}_m)$ be the expected demand. Algorithm \ref{alg:myopic_policy} describes the myopic policy.
\begin{algorithm}[h!]
	\caption{A myopic policy for the stochastic CSP}
	\label{alg:myopic_policy}
	\begin{algorithmic}[1]
		\Statex \textbf{input}: Expected demand vector $\bar{d}$, cutting patterns $a_j$, trim loss costs $g_j, \quad \forall j \in \{1,2, \dots, n\}$;
		\State Sample initial inventory $s_t \in \mathcal{S}$;
		\For{$t \gets 0 \dots$}
		\State Compute decision $x_t$ by solving the following integer linear  program:
		\begin{align}
		\text{min} \quad & \sum_{j=1}^n g_j x_{jt} \notag\\
		\text{s.t.} \quad & \notag \\
		& \sum_{j=1}^n a_{ij} x_{jt}+s_{it} \geq \bar{d}_i, \quad i \in \{1,2,\dots,m\}, \notag\\
		& x_{jt} \in \mathbb{Z}_+; \notag
		\end{align}
		\State Sample demand $d_{t+1} \sim p(d|s_t,x_t)$;
		\State Compute cost $c_{t+1} = c(s_t,x_t,d_{t+1})$ from equation  \eqref{eq:funcao_custo};
		\State Compute transition $s_{t+1} = f(s_t,x_t,d_{t+1})$ from equation \eqref{eq:funcao_transicao};
		\EndFor
	\end{algorithmic}
\end{algorithm}

\subsection{Dataset}
In the experiments, we use data originating from a real enterprise which faces the problem of cutting steel bars for building construction. Seven different lengths of steel bars (Table \ref{tab:item_types}) may be cut from stock bars with length 1500 cm. Demand for each bar length at each time period is stationary but random, so that the enterprise cannot predict exactly how much of each bar length customers will order. The enterprise wants to decide in anticipation which bar lengths to cut in order to meet demand without incurring high trim loss cost or maintaining high inventory.

The number of possible cutting patterns in this problem is large (all feasible solution to an integer knapsack problem with length 1500 cm and 7 item types). From a practical perspective, enterprises often restrict the number of cutting patterns to a small manageable set of \emph{good patterns} among efficient ones, i.e., those with small trim loss. Other selection criteria depend on  particular factors of the manufacturing process, such as  easiness of setup and ergonomics. We then have handcrafted 15 cutting patterns given in Table \ref{tab:padroes}, and as such they are not optimized to the enterprise operations. (Patterns used in practice were not disclosed by the enterprise.) 
\begin{table}
	\centering
	\caption{Lengths of demanded item types. Length of stock objects is 1500 cm}
	\begin{tabular}{cccccccc}
		\hline
		Item & 1 & 2 & 3 & 4 & 5 & 6 & 7 \\
		\hline
		Length (cm) & 115 & 180 & 267 & 314 & 880 & 1180 & 1200 \\
		\hline
	\end{tabular}
	\label{tab:item_types}
\end{table}
\begin{table}[t]
    \centering
    \caption{Cutting patterns used in the numerical experiments. Each column shows a different cutting pattern with associated trim losses}
    \begin{tabular}{cccccccccccccccc}
    \hline
         Item & P1 & P2 & P3 & P4 & P5 & P6 & P7 & P8 & P9 & P10 & P11 & P12 & P13 & P14 & P15\\
         \hline
         1 & 10 & 13 & 3 & 3 & 2 & 2 & 1 & 1 & 0 & 0 & 0 & 0 & 0 & 0 & 0 \\
         2 & 0  & 0  & 1 & 1 & 2 & 0 & 1 & 1 & 0 & 0 & 0 & 0 & 0 & 0 & 1 \\
         3 & 0  & 0  & 0 & 0 & 0 & 0 & 0 & 0 & 1 & 1 & 1 & 2 & 2 & 3 & 0 \\
         4 & 1  & 0  & 0 & 3 & 0 & 0 & 0 & 0 & 0 & 0 & 1 & 0 & 3 & 2 & 4 \\
         5 & 0  & 0  & 1 & 0 & 1 & 0 & 0 & 0 & 0 & 0 & 1 & 1 & 0 & 0 & 0 \\
         6 & 0  & 0  & 0 & 0 & 0 & 0 & 0 & 1 & 0 & 1 & 0 & 0 & 0 & 0 & 0 \\
         7 & 0  & 0  & 0 & 0 & 0 & 1 & 1 & 0 & 1 & 0 & 0 & 0 & 0 & 0 & 0 \\
         \hline
         Trim loss (cm) & 36 & 5 & 95 & 33 & 30 & 70 & 5 & 25 & 33 & 53 & 39 & 86 & 24 & 71 & 64 \\
         \hline
    \end{tabular}
    \label{tab:padroes}
\end{table}
We assumed the demand follows a mixture distribution in which the vector $d = (d_1, d_2, \dots, d_m)$ of demands for items follow a multinomial distribution $d \sim \mathtt{multinomial}(d_\text{total},p)$ conditional on the total demand $d_\text{total} \sim \mathtt{DiscUnif}(d_\text{min},d_\text{max})$, which follows a discrete uniform probability distribution, and $p = (p_1, p_2, \dots, p_m)$ is the vector of probabilities of each item being demanded. Thus, the marginal distribution of the vector $d$ is:
\begin{align}
\mathbb{P}(d) & = \sum_{d_\mathrm{total} = d_\text{min}}^{d_\mathrm{max}}  \mathbb{P}(d|d_\mathrm{total}) \mathbb{P}(d_\mathrm{total}), \notag \\
& = \sum_{d_\mathrm{total} = d_\text{min}}^{d_\mathrm{max}} \mathtt{Multinomial}(d_\text{total},p) \times \mathtt{DiscUnif}(d_\text{min},d_\text{max}). \label{eq:demand_probability}
\end{align}
\begin{table}[t]
	\centering
	\caption{Probability distribution of the demand for items}
	\begin{tabular}{cccccccc|c|c}
		\hline
		Item & 1 & 2 & 3 & 4 & 5 & 6 & 7 & $d_\text{min}$ & $d_\text{max}$  \\
		\hline
		Probability (p) & 0.30 & 0.20 & 0.20 & 0.10 & 0.10 & 0.05 & 0.05 & 40 & 50 \\
		\hline
	\end{tabular}
	\label{tab:demand}
\end{table}
\begin{table}[h!]
	\centering
	\caption{Additional problem data}
	\begin{tabular}{ccl}
		\hline
		Parameter & Value & Description \\
		\hline
		$h_i^+$ & $0.01 l_i$ & Unit cost of holding inventory, where $l_i$ is the length of item $i \in \{1,2,\dots, 7\}$\\
		$h_1^-$ & $1.0 l_i$ & Penalty unit cost of not meeting the demand\\
		$g_j$ & $0.1 g^{+}_j$ & Trim loss cost, in which $g^{+}_j$ is the trim loss of pattern $j \in \{1,2,\dots, 15\}$\\
		$s_\text{max}$ & 70 & Maximum inventory for each item in a time period \\
		$x_\text{max}$ & 30 & Number of stock objects available in each time period\\
		\hline
	\end{tabular}
\label{tab:problem_parameters}
\end{table}
\begin{table}[h!]
	\centering
	\caption{Parameters used in the numerical experiments}
	\begin{tabular}{ccl}
		\hline
		Parameter & Value & Description \\
		\hline
		$\gamma $& 0.8  & Discount factor \\
		$L_1$ & 30 & Number of policy iterations \\
		$L_2$ & $50 \times 10^3$ & Number of simulated state transitions \\
		$N_1$ & 10 & Number of iteration in the cross-entropy method \\
		$N_2$ & 100 &  Sample size of candidate solutions in the cross-entropy method \\
		\hline
	\end{tabular}
	\label{tab:experiment_parameters}
\end{table}
Samples from the probability distribution of the demand given in Eq. \eqref{eq:demand_probability} can be simulated simply by sampling the total demand $d_\text{total} \sim \mathtt{DiscUnif}(d_\text{min},d_\text{max})$ then sampling $d \sim \mathtt{Multinomial}(d_\text{total},p)$. Table \ref{tab:demand} shows the parameters of the probability function of the demand. Other parameters related to the problem are given in Table \ref{tab:problem_parameters}.

\subsection{Experiment results}
We implemented the algorithms in Python 3.7 with the aid of NumPy and Numba. The experiments were run in a Core i7 7700K machine with 8GB RAM. Each experiment run takes about 20 hours. In the experiments, we have used the parameters indicated in Table \ref{tab:experiment_parameters}. Initial values for $\theta_k^{(0)}, k \in {1, 2, \dots, K}$ are sampled from a Gaussian distribution with mean 0 and standard deviation 1. Notice that, since we set a sample size of $50 \times 10^3$ state transitions and 30 policy iterations, each experiment run takes a total of 1.5 million state transitions.
\begin{figure}
	\centering
	\includegraphics[scale=0.72]{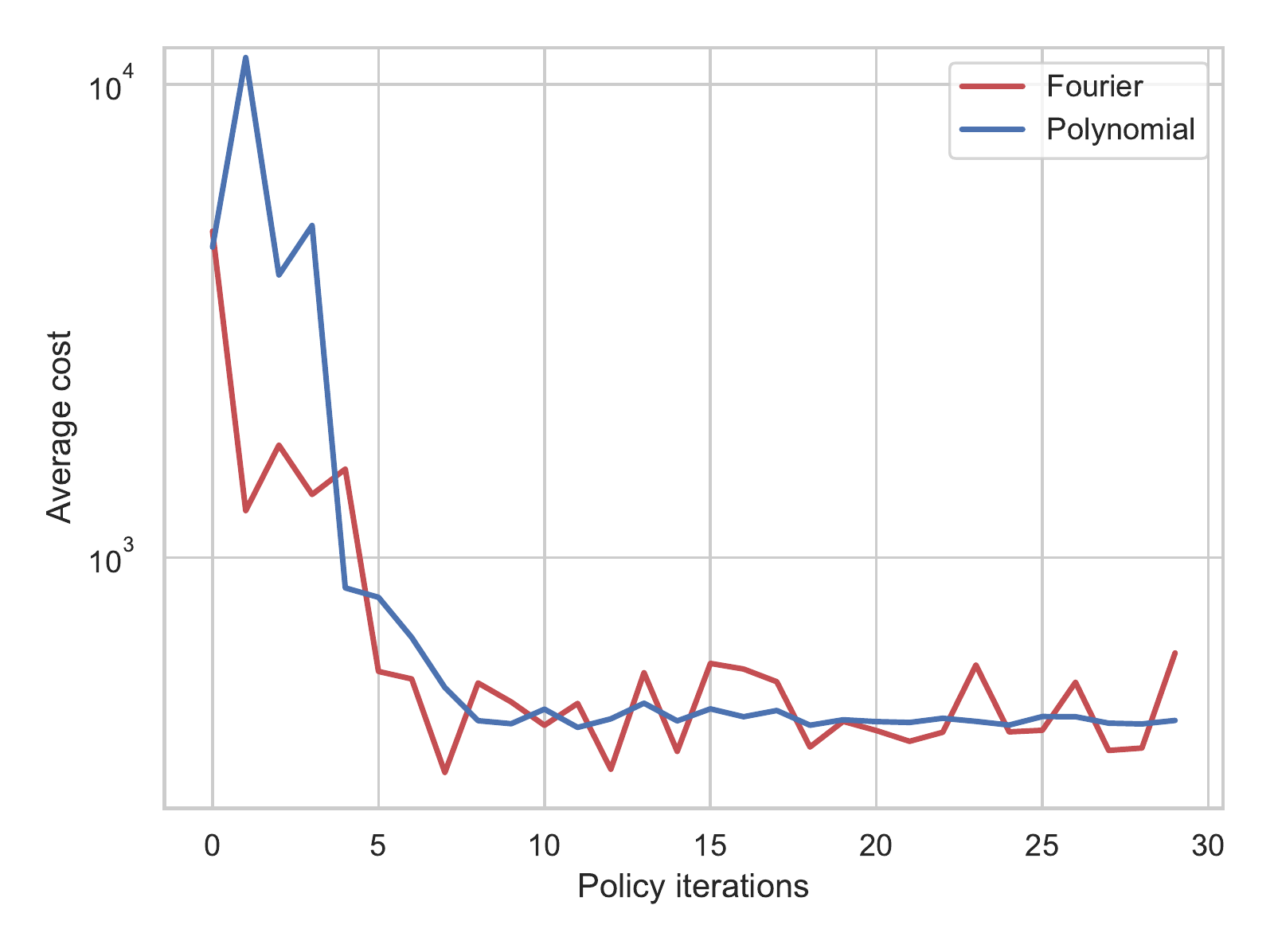}
	\caption{Average cost (log scale) of heuristic policies trained by our algorithm, evaluated with 10 simulation replications after training}
	\label{fig:training}
\end{figure}
\begin{figure}[h!]
	\centering
	\includegraphics[scale=0.63,trim={0 0.5cm 0 0.5cm},clip]{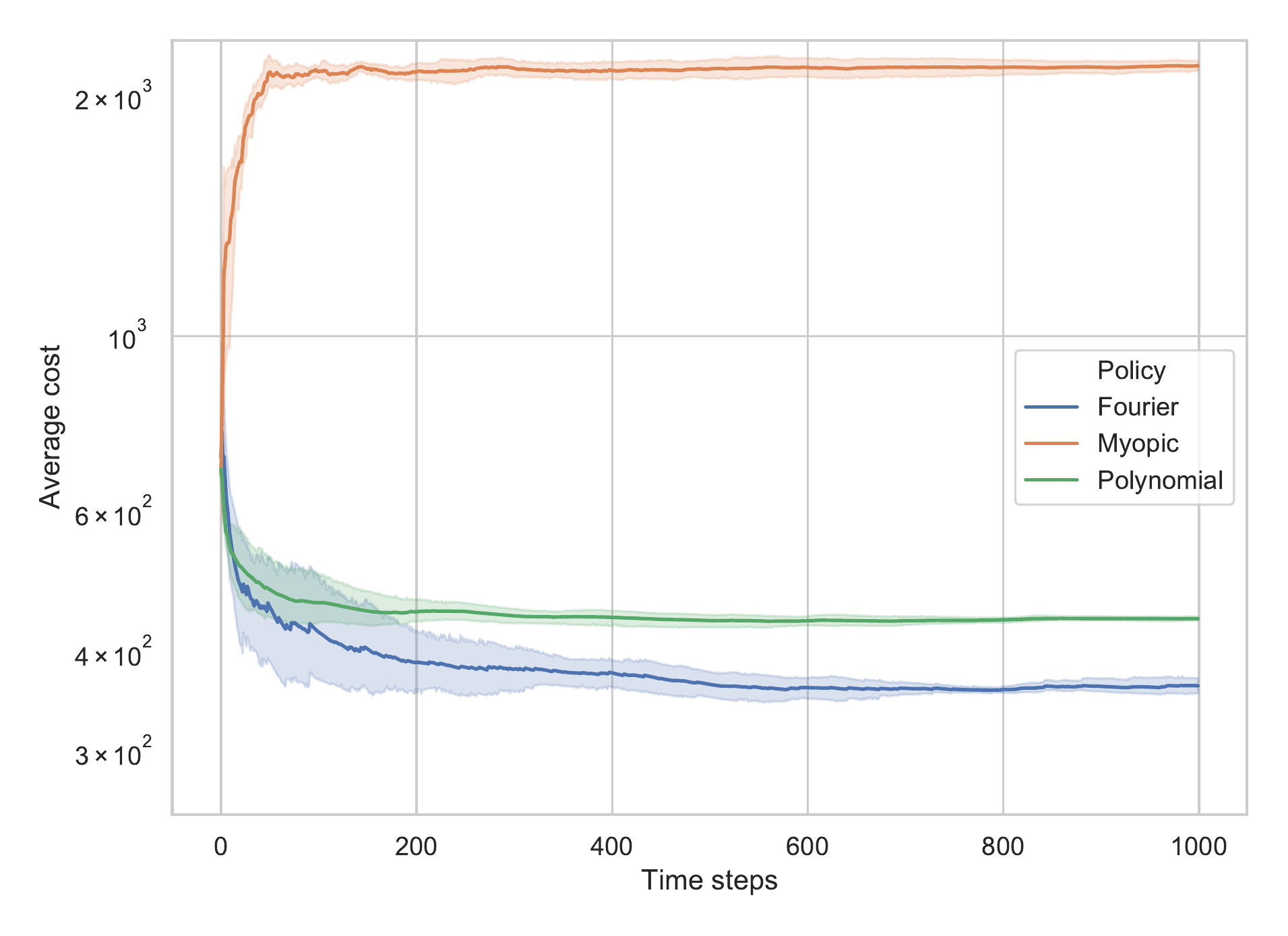}
	\caption{Average cost of policies (log scale). Bands are 95\% bootstraped confidence intervals for a sample of 10 simulation replications. Final average cost for Fourier and polynomial basis policies are 363.4 and 441.3, respectively, while for the myopic policy is 2186.5. The cost of the random policy is 6955.8 (omitted in the graph)} \label{fig:policies_rolling_cost}
\end{figure}

Due to simulation error and approximation of the action-value functions, the successive policies generated during training do not always exhibit monotonically better performance. Then, after training, we carry out a reevaluation of each of the 30 policies with 10 simulation replications in order to identify the best policy generated. Fig. \ref{fig:training} exhibits the performance of the successive policies after reevaluation, and it can be seen that policies improve in the first iterations and then start to oscillate. This policy oscillation phenomenon has been documented early in the literature from experimental studies in different applications of approximate policy iteration algorithms \citep{bertsekasNeuro}.

Trained policies were compared with the myopic policy described in Algorithm \ref{alg:myopic_policy} and a random policy, in which decisions are sampled uniformly in the set $\mathcal{X}_s$ at each time step. Fig. \ref{fig:policies_rolling_cost} shows the average cost of the policies over 10 replications of simulation of the demand. Trained policies achieved better performance than a myopic policy, with the policy trained using the Fourier basis showing lower average cost than the policy trained using the polynomial basis. Final average cost for Fourier and polynomial basis policies are 363.4 and 441.3, respectively, while for the myopic policy is 2186.5 and for the random policy is 6955.8. This represents an improvement in performance of approximately 80\% over the myopic policy.
\begin{figure}
	\centering
	\includegraphics[scale=0.7]{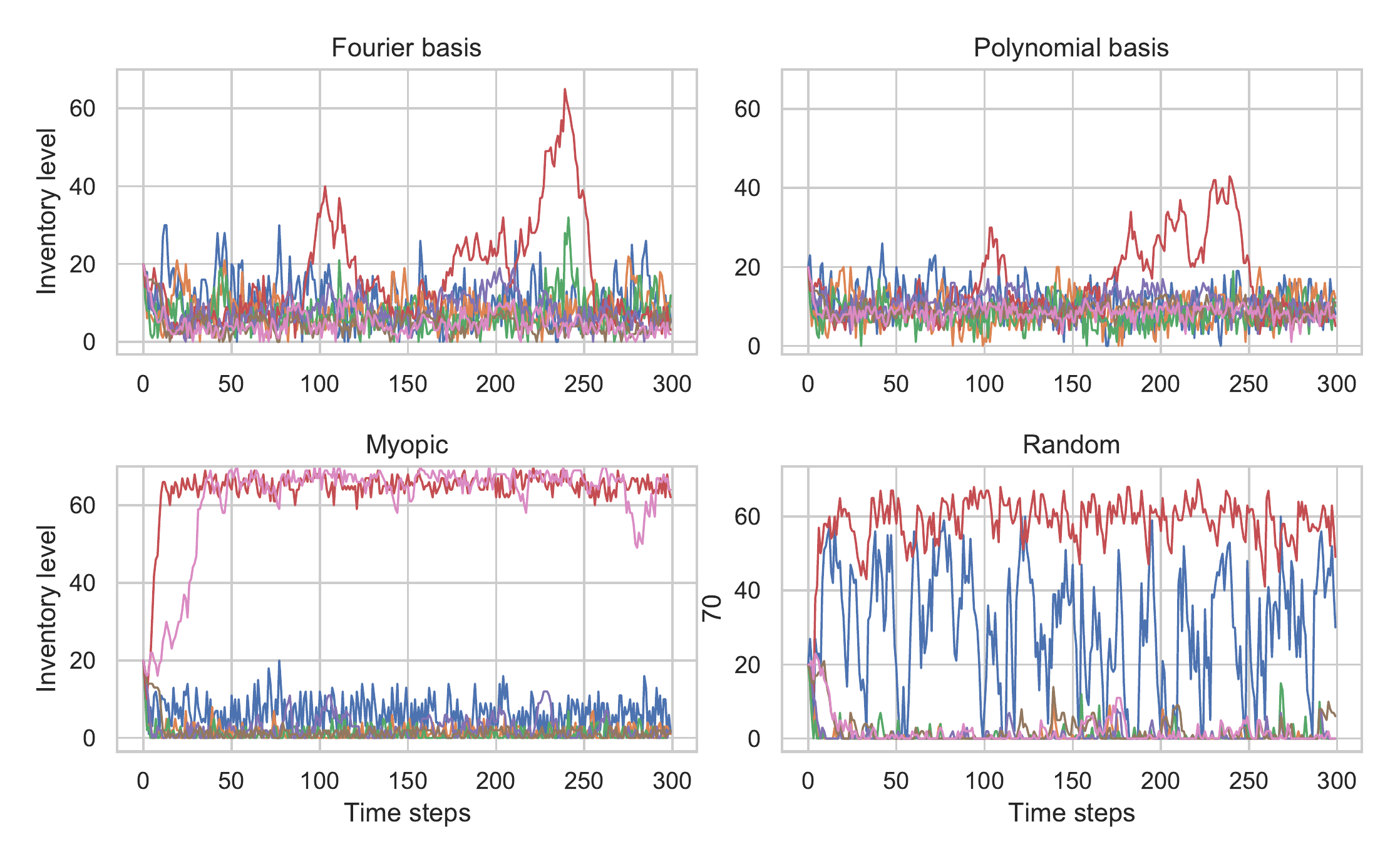}
	\caption{Initial inventory (state before decision) of all items at each time step. The myopic policy fails to adequately control inventory levels, while the trained policies maintain low inventory levels}
	\label{fig:estoques}
\end{figure}

Fig. \ref{fig:estoques} shows the initial inventory (state before decision) at each time step for each policy in a simulated trajectory of demands, states and decisions. Notice that the myopic policy fails to adequately control inventory levels, maintaining very low levels for some items and very high levels for others. We suppose this is due to the myopic nature of the policy, which does not take into account the impact of decisions in the future states. In contrast, the trained policies maintain sufficient levels to meet the demand while avoiding that the levels get high. However, as can be seen in Fig. \ref{fig:estoques}, sometimes inventory levels for some items may digress to high levels before the policy manage to get them back to low levels.

\begin{figure}
\centering
\includegraphics[scale=0.76,trim={0 0.5cm 0 0.5cm},clip]{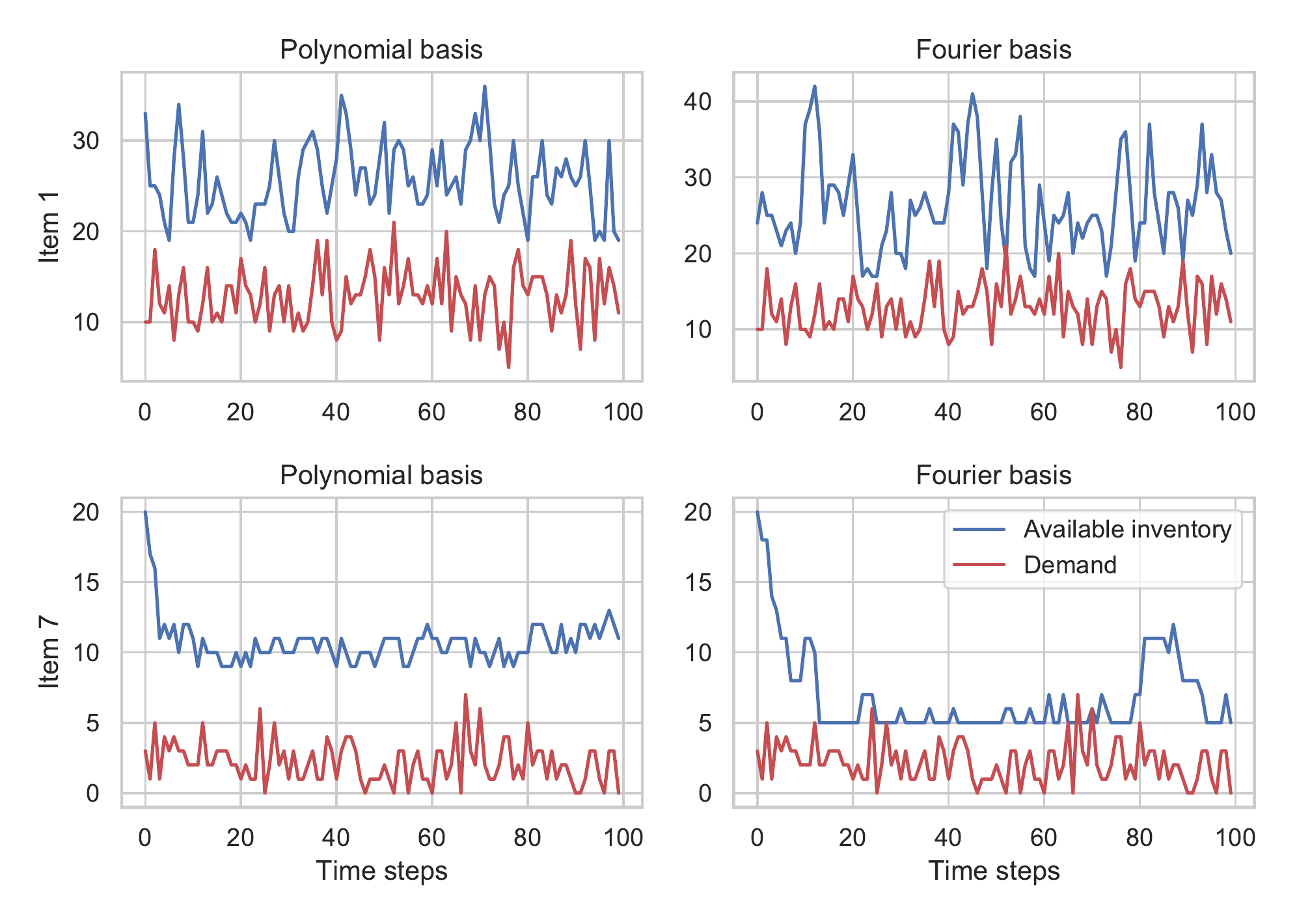}
\caption{Available inventory (sum of initial inventory and cut items) generated by the trained policies, and demand for items at each time step. It can be seen that the trained policies maintain sufficient available inventory in order to fulfill the demand}
\label{fig:available_inventory}
\end{figure}

Although both policies trained with a Fourier basis and a polynomial basis action-value function approximation exhibited lower average cost than a myopic policy, Fourier basis performed better. We can get some insight on this difference in performance when we look at the available inventory set by both policies over time. The available inventory is the post-decision state given by Eq. \eqref{eq:post-decision}, i.e., it is the sum of the initial inventory at a time period and the items produced by cutting stock objects. Fig. \ref{fig:available_inventory} shows the available inventory for items 1 and 7. Notice that the policy trained with the Fourier basis maintains lower available inventory while still sufficient to fulfill the demand. Interestingly, the policy with the Fourier basis function approximation  seems to have learned a \emph{safety} inventory level of 5 units for item 7, which has very low demand.
\begin{figure}[h!]
	\centering
	\includegraphics[scale=0.75,trim={0 0.5cm 0 0.5cm},clip]{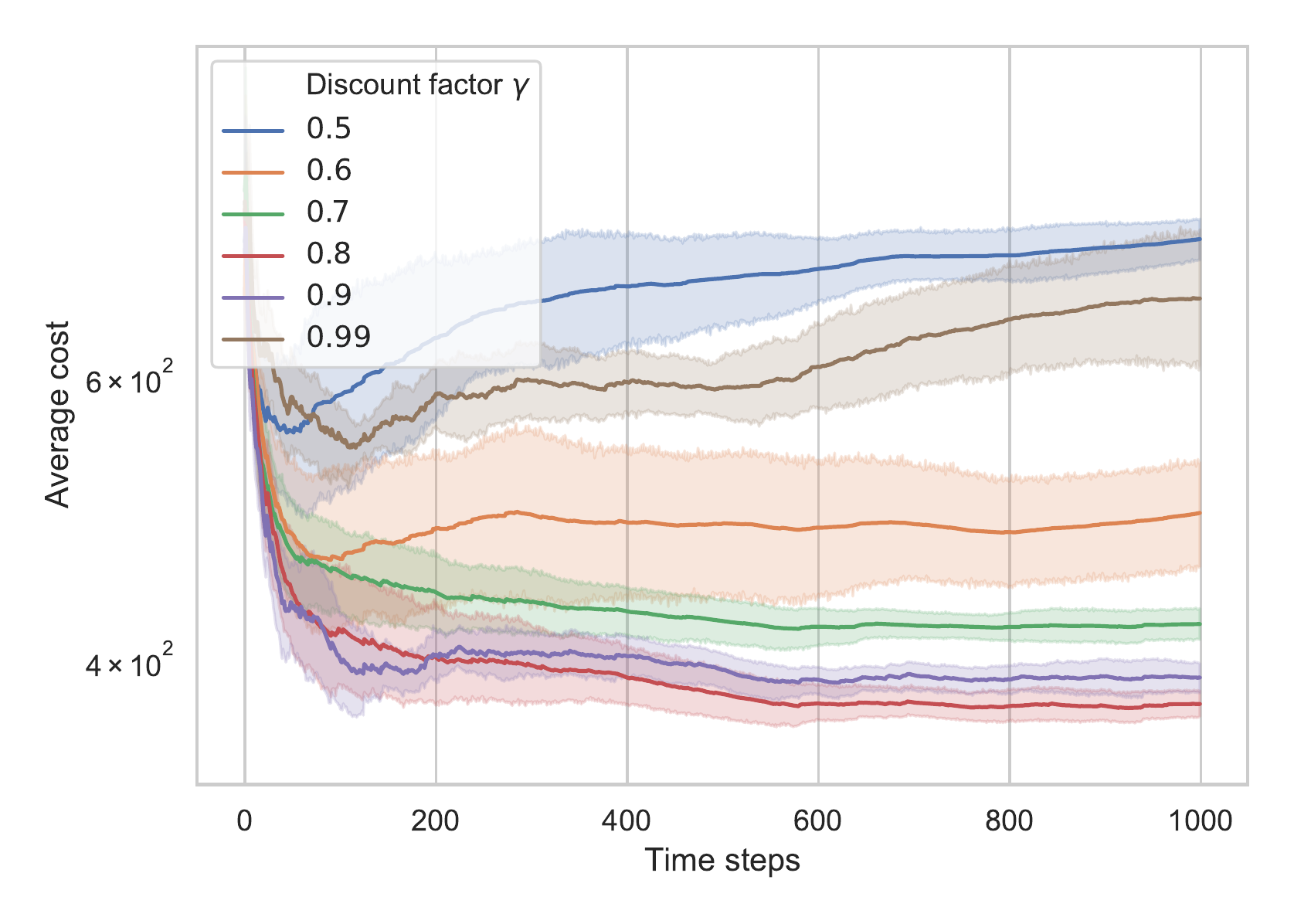}
	\caption{Performance of trained policies with the Fourier basis for different values of the discount factor $\gamma$. Performance is critically affected by the choice of discount factor}
	\label{fig:discount_factor}
\end{figure}

Finally, we have also investigated the effect of using different values for the discount factor $\gamma$, which works as a hyperparameter of the model. It is not obvious at first which value to assign to the discount factor and its effect on the performance of the trained policy. We have run experiments with the same dataset and 6 different values for the discount factor with a Fourier basis function approximation, with 10 simulation replications for each discount factor value. The results are given in Fig. \ref{fig:discount_factor}. It can be seen that the performance of the trained policy is very sensitive to different values of the discount factor. Lower levels (0.5 and 0.6) have considerably worse average cost than higher values, with lowest average cost at value 0.8. This indicates that there may be an optimal discount factor which results in the best performance. In practical applications, the best discount factor may be searched by using an algorithm for hyperparameter optimization.

\section{Conclusions}
In this paper, we have developed a solution approach to the stochastic cutting stock problem based on reinforcement learning. In the numerical experiments with realistic data, the obtained heuristic policies showed a considerable performance improvement of up to 80\% over a myopic policy which solves a deterministic formulation of the cutting stock problem. These promising results provide experimental evidence that a decision system for the cutting stock problem based on reinforcement learning may result in considerable cost reductions in industry. Nevertheless, more studies are necessary to evaluate the quality and stability of policies under different conditions, such as nonstationary demand and varying number of items or cutting patterns. Finally, we also see as possible research directions the investigation of nonlinear approximators such as artificial neural networks and the application of hyperparameter optimization.

{\small
\section*{Acknowledgements}
Funding: This study was financed in part by the Coordenação de Aperfeiçoamento de Pessoal de Nível Superior – Brasil (CAPES) – Finance Code 001; and Conselho Nacional de Desenvolvimento Científico e
Tecnológico (CNPq; Grant No.: 422464/2016-3). We also thank NVIDIA Corporation for the GPU support.}

\bibliography{references}

\end{document}